\documentclass[11pt]{amsart}

\usepackage{amsthm}
\usepackage{amsmath}
\usepackage{amssymb}

\usepackage[arrow,matrix,curve,cmtip,ps]{xy}

\newtheorem{thm}{Theorem}[section]   % Num
\newtheorem{cor}[thm]{Corollary}     % Numbered along with thm
\newtheorem{lem}[thm]{Lemma}         % Numbered along with thm
\newtheorem{prop}[thm]{Proposition}  % Numbered along with thm

\theoremstyle{definition}
\newtheorem{defn}[thm]{Definition}   % Numbered along with thm
 
\newtheorem{rem}[thm]{Remark}        % Numbered along with thm
\newtheorem{ex}[thm]{Example}        % Numbered along with thm

\numberwithin{equation}{section}

\newcommand{\Htil}{\tilde{H}}

\newcommand{\pvoh}{p_{v_0}^H}
\newcommand{\Prim}{\hbox{Prim\,}}

\begin{document}

\title[The $C^*$-algebras of arbitrary graphs]{The
$\boldsymbol{C^*}$-algebras of arbitrary graphs}

\author{D. Drinen}

\address{Dartmouth College\\Department of Mathematics\\Bradley Hall
6188\\Hanover, NH 03755}

\curraddr{Department of Mathematics \\ The University of the South\\
\\ Sewanee, TN 37383}

\email{ddrinen@sewanee.edu}

\author{M. Tomforde}

\address{Dartmouth College\\Department of Mathematics\\Bradley Hall
6188\\Hanover, NH 03755}

\curraddr{Department of Mathematics\\ University of Iowa\\
Iowa City\\ IA 52242\\ USA}

\email{tomforde@math.uiowa.edu}


\begin{abstract}

To an arbitrary directed graph we associate a row-finite directed graph whose
$C^*$-algebra contains the $C^*$-algebra of the original
graph as a full corner.  This allows us to generalize results for
$C^*$-algebras of row-finite graphs to $C^*$-algebras of arbitrary graphs: 
the uniqueness theorem, simplicity criteria, descriptions of the ideals and
primitive ideal space, and conditions under which a graph algebra is AF and
purely infinite.  Our proofs require only standard Cuntz-Krieger
techniques and do not rely on powerful constructs such as groupoids,
Exel-Laca algebras, or Cuntz-Pimsner algebras.

\end{abstract}

\maketitle

%------------------------------------------------------------------------
\section{Introduction}
\label{sec-intro}
%-----------------------------------------------------------------------

Since they were first introduced in 1947 \cite{Seg}, $C^*$-algebras 
have become important tools for mathematicians working in many areas.
Because of the immensity of the class of all $C^*$-algebras, however, it
has become important to identify and study special types of
$C^*$-algebras. These special types of $C^*$-algebras
(e.g.~AF-algebras, Bunce-Deddens algebras, AH-algebras, irrational
rotation algebras, group $C^*$-algebras, and various crossed products)
have provided great insight into the behavior of more
general $C^*$-algebras.  In fact, it is fair to say that much of the 
development of operator algebras in the last twenty years has been
based on a careful study of these special classes.

One important and very natural class of $C^*$-algebras comes
from considering $C^*$-algebras generated by partial isometries. 
There are a variety of ways to construct these $C^*$-algebras, but
typically any such construction will involve having the partial
isometries satisfy relations that describe how their initial and
final spaces are related.  Furthermore, one finds that in practice it
is convenient to have an object (e.g.~a matrix, a graph, etc.) that
summarizes these relations.

In 1977 Cuntz introduced a class of $C^*$-algebras that
became known as Cuntz algebras \cite{Cun}.  For each $n=2,3,\ldots,
\infty$ the Cuntz algebra $\mathcal{O}_n$ is generated by $n$
isometries satisfying certain relations.  The Cuntz algebras were
important in the development of $C^*$-algebras because they provided
the first examples of $C^*$-algebras whose $K$-theory has torsion.  In
1980 Cuntz and Krieger considered generalized versions of the Cuntz
algebras \cite{CK1}.  Given an $n \times n$ matrix $A$ with
entries in $\{0,1\}$, the Cuntz-Krieger algebra $\mathcal{O}_A$ is
defined to be the $C^*$-algebra generated by partial isometries
satisfying relations determined by $A$.  A study of the Cuntz-Krieger
algebras was made in the seminal paper \cite{CK1} where it was shown
that they arise naturally in the study of topological Markov chains. 
It was also shown that there are important parallels between
these $C^*$-algebras and certain kinds of dynamical systems
(e.g.~shifts of finite type).

In 1982 Watatani noticed that by considering a $\{0,1\}$-matrix as
the adjacency matrix of a directed graph, one could view Cuntz-Krieger
algebras as $C^*$-algebras associated to certain finite directed
graphs \cite{Wat}.  Although Watatani published some papers using
this graph approach \cite{FW,Wat}, his work went largely unnoticed. 
It was not until 1997 that Kumjian, Pask, Raeburn, and Renault
rediscovered $C^*$-algebras associated to directed graphs.  

This theory of $C^*$-algebras associated to graphs was developed in
\cite{kprr}, \cite{kpr}, and \cite{bprs}.  In these papers the authors
were able to define and work with $C^*$-algebras associated
to finite graphs as well as $C^*$-algebras associated to infinite
graphs that are row-finite (i.e.~all vertices emit a finite number of
edges).  By allowing all finite graphs as well as certain infinite
graphs, these graph algebras included many $C^*$-algebras that were
not Cuntz-Krieger algebras. Furthermore, it was found that the graph
not only described the relations for the generators, but also many
important properties of the associated $C^*$-algebra could be
translated into graph properties.  Thus the graph provides a tool for
visualizing many aspects of the associated $C^*$-algebra.  In
addition, because graph algebras consist of a wide class of
$C^*$-algebras whose structure can be understood, other areas of
$C^*$-algebra theory have benefitted nontrivially from their study.

Despite these successes, many people were still unsatisfied with the
condition of row-finiteness and wanted a theory of
$C^*$-algebras for arbitrary graphs.  This desire was
further fueled by the fact that in his original paper
\cite{Cun} Cuntz defined a $C^*$-algebra $\mathcal{O}_\infty$,
which seemed as though it should be the $C^*$-algebra associated to a
graph with one vertex and a countably infinite number of edges. 
Despite many people's desire to extend the definition of
graph algebras to arbitrary graphs, it was unclear exactly
how to make sense of the defining relations in the non-row-finite
case.  It was not until 2000 that Fowler, Laca, and
Raeburn were finally able to extend the definition of graph algebras
to arbitrary directed graphs \cite{flr}.  These graph algebras now
included the Cuntz algebra $\mathcal{O}_\infty$, and as
expected it arises as the $C^*$-algebra of the graph with
one vertex and infinitely many edges.

In the time since $C^*$-algebras associated to
arbitrary graphs were defined, there have been many
attempts to extend results for row-finite graph algebras to
arbitrary graph algebras.  However, because many of the proofs of the
fundamental theorems for $C^*$-algebras of row-finite graphs make
heavy use of the row-finiteness assumption, it has often been unclear
how to proceed. In most cases where results have been generalized,
the proofs have relied upon sophisticated techniques and powerful
machinery such as groupoids, the Exel-Laca algebras of \cite{el}, and
the Cuntz-Pimsner algebras of \cite{pimsner}.

In this paper we describe an operation called desingularization that
transforms an arbitrary graph into a row-finite graph with no
sinks.  It turns out that this operation preserves Morita equivalence
of the associated $C^*$-algebra as well as the loop structure and path
space of the graph.   Consequently, it is a powerful tool in the
analysis of graph algebras because it allows one to apply much of the
machinery that has been developed for row-finite graph algebras to
arbitrary graph algebras.

Desingularization was motivated by the process of ``adding a
tail to a sink" that is described in \cite{bprs}.  In fact,
this process is actually a special case of desingularization. 
The difference is that now we not only add tails at sinks, but
we also add (more complicated) tails at vertices that emit
infinitely many edges.  Consequently, we shall see that vertices
that emit infinitely many edges will often behave similarly to
sinks in the way that they affect the associated $C^*$-algebra.  In
fact for some of our results, such as conditions for simplicity, one
can take the result for row-finite graphs and replace the word
``sink" by the phrase ``sink or vertex that emits infinitely many
edges" to get the corresponding result for arbitrary graphs.

We begin in Section \ref{sec-desing} with the definition of
desingularization.  This is our main tool for dealing with
$C^*$-algebras associated to arbitrary graphs.  It gives the reader
who is comfortable with $C^*$-algebras of row-finite graphs a great
deal of intuition into the structure of non-row-finite graph
algebras.  This is accomplished by providing a method for easily
translating questions about arbitrary graph algebras to the
row-finite setting. After the definition of desingularization, we
describe a correspondence between paths in the original graph and
paths in the desingularization.  We then show that desingularization
preserves loop structure of the graph as well as Morita equivalence
of the $C^*$-algebra.  This allows us to obtain easy proofs of
several known results.  In particular, we prove the uniqueness
theorem of \cite{flr} and give necessary and sufficient conditions for
a graph algebra to be simple, purely infinite, and AF.

In Section~\ref{sec-ideal} we describe the ideal structure of graph
algebras.  Here we will see that our solution is more complicated than
what occurs in the the row-finite case.  The correspondence with
saturated hereditary sets described in \cite{bprs} no longer holds. 
Instead we have a correspondence of the ideals with pairs $(H,S)$,
where $H$ is a saturated hereditary set and $S$ is a set containing
vertices that emit infinitely many edges, only finitely many of which
have range outside of $H$.  

We conclude in Section~\ref{sec-prim} with a description of the
primitive ideal space of a graph algebra.  Our result will again be
more complicated than the corresponding result for the row-finite
case, which involves maximal tails \cite{bprs}.  For arbitrary graphs
we will need to account for vertices that emit infinitely many edges,
and our description of the primitive ideal space will include both
maximal tails and special vertices that emit infinitely many edges
known as ``breaking vertices".

We thank Iain Raeburn for making us aware of the related papers
by Szyma\'nski \cite{sz} and Paterson \cite{patgg}, and we thank both
Iain Raeburn and Dana Williams for their comments on the first draft
of this paper.  After this work was completed, it was brought to our
attention that our description of the primitive ideal space in
Section~\ref{sec-prim} had been obtained independently in the preprint
\cite{bhrs}.  Although the results in \cite{bhrs} are similar to
some of our results in Section~\ref{sec-prim}, one should note that
the methods used in the proofs are very different.  In addition, we
mention that we have adopted their term ``breaking vertex" to provide
consistency for readers who look at both papers.

%------------------------------------------------------------------------
\section{The desingularized graph}
\label{sec-desing}
%-----------------------------------------------------------------------

We closely follow the notation established in \cite{kpr} and
\cite{bprs}.  A (directed) graph $E=(E^0, E^1, r, s)$ consists of
countable sets $E^0$ of vertices and $E^1$ of edges, and maps $r,s:E^1
\rightarrow E^0$ describing the source and range of each edge.  
We let $E^*$ denote the set of finite paths in $E$, and we let $E^\infty$
denote the set of infinite paths.  The maps $r,s$ extend to $E^*$ in the
obvious way and $s$ extends to $E^\infty$.  

A vertex $v$ is called a \emph{sink} if $|s^{-1}(v)|=0$, and $v$ is called an
\emph{infinite-emitter} if $|s^{-1}(v)|=\infty$.  If $v$ is either a sink or
an infinite-emitter, we call it a {\it singular vertex}.  A graph $E$ is
said to be {\it row-finite} if it has no infinite-emitters.

Given any graph (not necessarily row-finite), a {\it
Cuntz-Krieger $E$-family} consists
of mutually orthogonal projections $\{p_v \,|\, v \in E^0\}$ and partial
isometries
$\{s_e \,|\, e \in E^1\}$ with orthogonal ranges satisfying the {\it
Cuntz-Krieger relations}:
\begin{enumerate}
\item $s_e^* s_e = p_{r(e)}$ for every $e \in E^1$;
\item $s_e s_e^* \leq p_{s(e)}$ for every $e \in E^1$;
\item $p_v = \sum_{\{e \,|\, s(e)=v\}} s_e s_e^*$ for every $v \in E^0$
that is not a singular vertex.
\end{enumerate}
The graph algebra $C^*(E)$ is defined to be the $C^*$-algebra generated
by a universal Cuntz-Krieger $E$-family.  For the existence of such a
$C^*$-algebra, one can either modify the proofs in
\cite[Theorem 2.1]{aHr} or \cite[Theorem 1.2]{kpr}, or one can appeal to
more general constructions such as \cite{blackadar} or
\cite{pimsner}.

Given a graph $E$ we shall construct a graph $F$, called a {\it
desingularization of $E$}, with the property that $F$ has no singular vertices
and $C^*(E)$ is isomorphic to a full corner of $C^*(F)$.  Loosely
speaking, we will build $F$ from $E$ by replacing every singular vertex $v_0$
in $E$ with its own infinite path, and then redistributing the edges of
$s^{-1}(v_0)$ along the vertices of the infinite path.  Note that if
$v_0$ happens to be a sink, then $|s^{-1}(v_0)| = 0$ and there are no edges
to redistribute.  In that case our procedure will coincide with the process of
adding an infinite tail to a sink described in \cite[(1.2)]{bprs}.  

\begin{defn}
\label{defn-addtail}
Let $E$ be a graph with a singular vertex $v_0$.  We {\it add a tail}
to $v_0$ by performing the following procedure.  If $v_0$ is a sink, we add a
graph of the form 
\begin{equation}
\label{tail}
\xymatrix{
v_0 \ar[r]^{e_1} & v_1 \ar[r]^{e_2} & v_2 \ar[r]^{e_3} & v_3 \ar[r]^{e_4}
& \cdots\\
}
\end{equation}
as described in \cite[(1.2)]{bprs}.  If $v_0$ is an infinite
emitter we first list the edges $g_1, g_2, g_3, \ldots$ of $s^{-1}(v_0)$.  Then
we add a graph of the form shown in (\ref{tail}), remove the edges in
$s^{-1}(v_0)$, and for every $g_j \in s^{-1}(v_0)$ we draw an edge $f_j$ from
$v_{j-1}$ to $r(g_j)$.

For any $j$ we shall also define $\alpha^j$ to be the path $\alpha^j := e_1e_2
\ldots e_{j-1}f_j$ in $F$.
\end{defn}

Note that different orderings of the edges of $s^{-1}(v_0)$ may give rise
to nonisomorphic graphs via the above procedure.

\begin{defn}
If $E$ is a directed graph, a {\it desingularization of $E$} is a
graph $F$ obtained by adding a tail at every singular vertex of $E$. 
\end{defn}

\begin{ex} 
\label{ex-mainex}
Suppose we have a graph $E$ containing this fragment:

$$
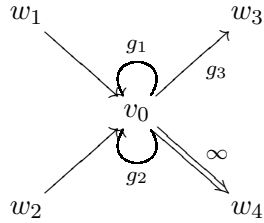
\xymatrix{
w_1 \ar[dr] & & w_3 \\
& v_0 \ar@{=>}[dr]^<>(.6){\infty} \ar[ur]_<>(.6){g_3}
\ar@(ur,ul)[]_{g_1}
\ar@(dr,dl)[]^{g_2} &\\
w_2 \ar[ur] & & w_4\\
}
$$

\vspace*{.1in}

\noindent where the double arrow labeled $\infty$ denotes a countably infinite
number of edges from $v_0$ to $w_4$.  Let us label the edges from $v_0$ to
$w_4$ as
$\{g_4, g_5, g_6, \ldots\}$.  Then a desingularization of $E$ is given by the
following graph $F$.

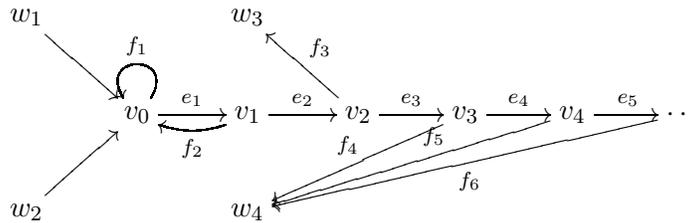
$$
\xymatrix{
w_1 \ar[dr] & & w_3 & & &\\
 & v_0 \ar@(ur,ul)[]_{f_1} \ar[r]^{e_1} & v_1 \ar[r]^{e_2}
\ar@/^/[l]^{f_2} & v_2 \ar[r]^{e_3}
\ar[ul]_{f_3} & v_3 \ar[r]^{e_4} \ar[dll]_{f_4} & v_4 \ar[r]^{e_5}
\ar[dlll]_<>(.4){f_5} & \cdots \ar[dllll]^{f_6}\\ 
w_2 \ar[ur] & & w_4 & & & &\\
}
$$
\end{ex}

\begin{ex}
\label{ex-Oinfty}
If $E$ is the $\mathcal{O}_{\infty}$ graph (one vertex with infinitely 
many loops), a desingularization $F$ looks like this:
$$
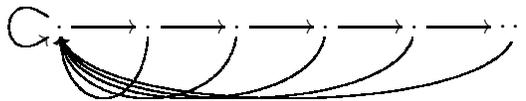
\xymatrix{
. \ar[r] \ar@(ul,dl)[] & . \ar[r] \ar@(d,d)[l] & . \ar[r] \ar@(d,d)[ll] &
. \ar[r] 
\ar@(d,d)[lll] & . \ar[r] \ar@(d,d)[llll] & \ar@(d,d)[lllll] \cdots\\
&&&&&\\
}
$$
\end{ex}

\begin{ex}
\label{ex-K}
The following graph was mentioned in \cite[Remark 11]{flr}:
$$
\xymatrix{
\cdots \ar[r] & . \ar[r] & . \ar[r] & v_0 \ar@{=>}[r]^{\infty}
& . \ar[r] & . \ar[r] & \cdots \\
}
$$
A desingularization of it is:
$$
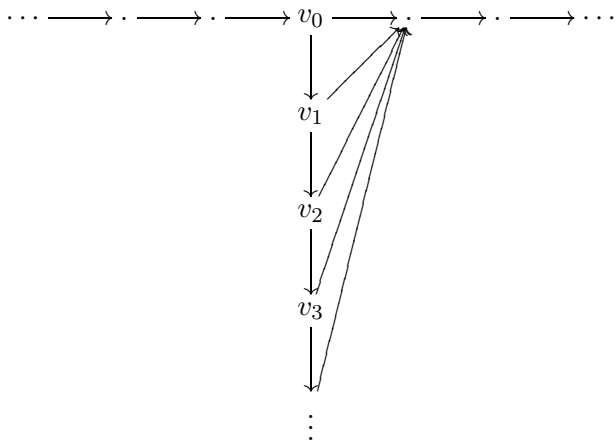
\xymatrix{
\cdots \ar[r] & . \ar[r] & . \ar[r] & v_0 \ar[r] \ar[d] & . 
\ar[r] & .
\ar[r] & \cdots\\
&&& v_1 \ar[d] \ar[ur] &&&\\
&&& v_2 \ar[d] \ar[uur] &&&\\
&&& v_3 \ar[d] \ar[uuur] &&&\\
&&& \vdots \ar[uuuur] &&&\\
}
$$
\end{ex}

It is crucial that desingularizing a graph preserves connectivity, path
space, and loop structure in the appropriate senses, and this will turn
out to be the case. We make these
ideas precise with the next three lemmas:  Lemma~\ref{lem-correspondence}
describes how the path spaces of $E$ and $F$ are related, Lemma~\ref{lem-LK}
shows that desingularization
preserves loop structure, and Lemma~\ref{lem-cofinalplus} describes the
relationship between cofinality of a graph and cofinality of its
desingularization.

We first review some notation.  If $E$ is a directed graph and $S_1,S_2
\subseteq E^0$ we say {\it $S_1$ connects to $S_2$}, denoted $S_1 \geq
S_2$, if for every $v \in S_1$ there exists $w \in S_2$ and $\alpha \in
E^*$ with $s(\alpha)=v$ and $r(\alpha)=w$.  Frequently one or both of
the $S_i$'s will contain a single vertex $v$, in which case we write $v$
rather than $\{v\}$.  If $\lambda$ is a finite or infinite path in $E$, we
write $S \geq \lambda$ to mean $S \geq \{s(\lambda_i)\}_{i=1}^{|\lambda|}$.
Finally, a graph $E$ is said to be {\it cofinal} if for every infinite path
$\lambda$ we have $E^0 \geq \lambda$.

\begin{lem}
\label{lem-correspondence}

Let $E$ be a graph and let $F$ be a desingularization of $E$.
\begin{itemize}
\item [(a)] There are bijective maps
$$
\phi: E^* \longrightarrow \{\beta \in F^* \,|\, s(\beta),
r(\beta) \in E^0\}\\ 
$$
$$
\phi_\infty: E^\infty \cup \{\alpha \in E^* \,|\, r(\alpha) \hbox{ is a
singular
vertex}\} \longrightarrow \{\lambda \in F^\infty \,|\,
s(\lambda) \in E^0\}.\\ 
$$
The map $\phi$ preserves source and range (and hence $\phi$ preserves loops),
and the map $\phi_\infty$ preserves source.
\item [(b)] The map $\phi_\infty$ preserves $\geq$ in the following sense.
For every $v \in E^0$ and $\lambda \in E^\infty \cup  \{ \alpha \in E^* \, |
\, r(\alpha) \text{ is a singular vertex} \}$, we have $v \geq \lambda$ in
$E$ if and only if $v \geq \phi_\infty(\lambda)$ in $F$. 
\end{itemize}
\end{lem}

\begin{proof}
First define a map $\phi':E^1 \rightarrow F^*$.  If $e \in
E^1$, then $e$ will have one of two forms:  either $s(e)$ is not
a singular vertex,
in which case $e \in F^1$, or else $s(e)$ is a singular vertex, in which
case
$e=g_j$ for some $j$.  We define $\phi'$ by 
$$
\phi'(e) = \left\{\begin{array}{ll}
e & \hbox{if $s(e)$ is not singular;}\\
\alpha^j & \hbox{if $e=g_j$ for some $j$,}\\
\end{array}\right.
$$
where $\alpha^j := e_1 \ldots e_{j-1}f_j$ is the path described in
Definition~\ref{defn-addtail}.  Since $\phi'$ preserves source and range, it
extends to a map on the finite path space $E^*$.  In particular, for $\alpha
 = \alpha_1 \ldots \alpha_n \in E^*$ define
$\phi(\alpha) =
\phi'(\alpha_1)\phi'(\alpha_2)\ldots \phi'(\alpha_{|\alpha|})$.  It is easy
to check that $\phi$ is
injective, that it preserves source and range, and that it is onto the
set $\{\beta \in F^* \,|\, s(\beta), r(\beta) \in E^0\}$.  We
define $\phi_\infty$
similarly. In particular, if $\lambda = \lambda_1 \lambda_2 \ldots \in
E^\infty$, define $\phi_\infty(\lambda) = \phi'(\lambda_1)\phi'(\lambda_2)
\ldots$. If $\alpha$ is a finite path whose range is a singular vertex $v_0$,
we define $\phi_\infty(\alpha) = \phi(\alpha)e_1e_2\ldots$, where
$e_1e_2,\ldots$ is the tail in $F$ added to $v_0$.  

To show that $\phi_\infty$ is a bijection, we construct an inverse
$\psi_\infty:\{\lambda \in F^\infty \,|\, s(\lambda) \in E^0\}
\rightarrow E^\infty \cup \{\alpha \in E^* \,|\,
r(\alpha) \hbox{ is a singular vertex}\}$.  Notice that every
$\lambda \in F^\infty$ either returns to $E$ infinitely often or it ends
up in one of the added infinite tails.  More precisely, 
$\lambda$ has one of two forms:  either $\lambda =
a^1a^2\ldots$ or $\lambda = a^1a^2 \ldots a^n e_1 e_2 e_3 \ldots$,   
where each $a^k$ is either an edge of $E$ or an $\alpha^j$. 
We define $\psi'$ by 
$$
\psi'(a^k) := \left\{\begin{array}{ll}
a^k & \hbox{if $a^k \in E^1$;}\\
g_j & \hbox{if $a^k=\alpha^j$ for some $j$,}\\
\end{array}\right.
$$
and we define
$$
\psi_\infty(\lambda) := \left\{\begin{array}{ll}
\psi'(a^1)\psi'(a^2)\ldots & \hbox{if $\lambda = a^1a^2\ldots$;}\\
\psi'(a^1)\ldots\psi'(a^n) & \hbox{if $\lambda = a^1 \ldots
a^ne_1e_2\ldots$.}\\
\end{array}\right.
$$
It is easy to check that $\phi_\infty$ and $\psi_\infty$ are inverses,
and we have established (a).

To prove (b), let $\lambda \in E^\infty \cup \{ \alpha \in E^* \,|\, r(\alpha)
\hbox{ is a singular vertex}\}$ and $v \geq \lambda$ in $E$.  Then there exists
a finite path $\alpha$ in $E$ such that $s(\alpha)=v$ and $r(\alpha)=w$
for some $w \in E^0$ lying on the path $\lambda$.  Note that the
vertices of $E$ that are on the path $\phi_\infty(\lambda)$ are
exactly the same as the vertices on the path $\lambda$.  Hence $w$ must
also be a vertex on the path $\phi_\infty(\lambda)$.  Now, because
$\phi$ preserves source and range, $\phi(\alpha)$ is a path that
starts at $v$ and ends at $w$, which is a vertex on $\phi_\infty(\lambda)$.  
Thus $v \geq \phi_\infty(\lambda)$.

For the converse let $\lambda \in E^\infty \cup \{ \alpha \in E^* \,|\, r(\alpha)
\hbox{ is a singular vertex}\}$ and $v \in E^0 $, and suppose that $v \geq
\phi_\infty(\lambda)$ in $F$.  Then there exists a finite path $\beta$ in $F$
with $s(\beta)=v$
and $r(\beta) = w$ for some vertex $w$ on the path
$\phi_\infty(\lambda)$.  Notice that if $r(\beta)$ is a vertex on one
of the added infinite tails, then $\phi_\infty(\lambda)$ must have passed
through $v_0$, and so must have $\beta$.  Thus we may assume $r(\beta) \in E^0
\subseteq F^0$.  Now $\beta$ is a finite path in $F$ that starts and ends in
$E^0$, so it can be pulled back to a path $\phi^{-1}(\beta) \in E^*$
with source $v$ and range $r(\beta)$.  Since $r(\beta)$ lies on the
path $\phi_\infty(\lambda)$, it lies on the path $\lambda$, and thus
$\phi^{-1}(\beta)$ is a path from $v$ to some vertex of $\lambda$.  Hence
$v \geq \lambda$ in $E$.
\end{proof}

A {\it loop} in a graph $E$ is a finite path $\alpha = \alpha_1 \alpha_2
\ldots \alpha_{|\alpha|}$ with $s(\alpha)=r(\alpha)$.  The vertex
$s(\alpha)=r(\alpha)$ is called the {\it base point} of the loop. A loop
is said to be {\it simple} if $s(\alpha_i)=s(\alpha_1)$ implies
$i=1$.  Therefore a simple loop is one that does not return to its base
point more than once.  An {\it exit} for a loop $\alpha$ is an edge $f$
such that $s(f)=s(\alpha_i)$ for some $i$, and $f \neq \alpha_i$.  A
graph $E$ is said to satisfy {\it Condition~(L)} if every loop has an
exit and $E$ is said to satisfy {\it Condition~(K)} if no vertex in $E$
is the base point of exactly one simple loop.

\begin{lem}
\label{lem-LK}

Let $E$ be a graph and let $F$ be a desingularization of $E$.  Then 
\begin{itemize} 
\item[(a)] $E$ satisfies Condition~(L) if and only if $F$ satisfies
Condition~(L).  
\item[(b)] $E$ satisfies Condition~(K) if and
only if $F$ satisfies Condition~(K).
\end{itemize}
\end{lem}

\begin{proof}
If $\alpha$ is a loop in $E$ with no exits, then all the vertices on
$\alpha$ emit exactly one edge.  Hence none of these vertices are singular
vertices, and $\phi(\alpha)$ is a loop
in $F$ with no exits. If $\alpha$ is a loop in $F$ with
no exits, then we claim that none of the singular vertices of $E$ can
appear in the loop.  
To see this, note that if $v_0$ is a sink in $E$, then it cannot be a part of
a loop in $F$; and if $v_0$ is an infinite-emitter in $E$, then $v_0$ is the
source of two edges, which would necessarily create an exit for any loop.
Since none of the singular vertices of $E$ appear in $\alpha$, it
follows that $\phi^{-1}(\alpha)$ is a loop in $E$ with no exits. This
establishes part (a).

Now suppose $v \in E^0$ is the base of exactly one simple loop $\alpha$ in
$E$. Then $\phi(\alpha)$ is a simple loop in $F$.  If there were another
simple loop $\beta$ in $F$ based at $v$, then $\phi^{-1}(\beta)$ would be
simple loop in $E$ based at $v$ that is different from $\alpha$.  Thus 
if $F$ satisfies Condition~(K), then $E$ satisfies Condition~(K). 

Now suppose $E$ satisfies Condition~(K).  Let $v \in F^0$ be the base of
a simple loop $\alpha$ in $F$.  If $v \in E^0$, then $\phi^{-1}(\alpha)$
is a simple loop in $E$ based at $v$.  Since $E$ satisfies Condition~(K),
there is a simple loop $\beta$ in $E$ different from $\phi^{-1}(\alpha)$.
Certainly, $\phi(\beta)$ is a simple loop in $F$ and, because $\phi$ is
injective, $\phi(\beta)$ must be different from $\alpha$.  

Now suppose $v$ is on an added infinite tail; that is, $v=v_n$ for some $n
\geq 1$. Then $\alpha$ must have the form $\alpha' e_1 e_2 \ldots e_n$ for
some $\alpha' \in F^*$. Now, $e_1 e_2 \ldots e_n \alpha'$ is a simple loop in
$F$ based at $v_0$ and hence $\phi^{-1}(e_1 \ldots e_n
\alpha')$ is a simple loop in $E$ based at $v_0$.  Since $E$ satisfies
Condition~(K), there must be another simple loop $\beta$ in $E$ based at $v_0$. 
Now $\phi(\beta)$
will be a simple loop in $F$ based at $v_0$.
If $v_n$ is not a vertex on $\phi(\beta)$, then $\alpha' \phi(\beta) e_1
\ldots e_n$ will be another simple loop based at $v_n$ that is different
from $\alpha$.
On the other hand, if $v_n$ is a vertex of $\phi(\beta)$, then $\phi(\beta)$
has the form $e_1 \ldots e_n \beta'$, where $\beta' \in F^*$.  Since
$\phi(\beta)$ is a simple loop based at $v_0$, we know that $s(\beta_i) \neq
v_0$ for $1 \leq i \leq |\beta'|$.  Hence $v_n$ is not a vertex on the path
$\beta'$.  Therefore $\beta'e_1 \ldots e_n$ is a simple loop based at $v_n$.
Furthermore, it is different from the loop $\alpha = \alpha' e_1 \ldots e_n$,
because if they were equal then we would have $\alpha' = \beta'$, which
contradicts the fact that $\alpha = \alpha'e_1 \ldots e_n$ and $\phi(\beta) =
\beta'e_1 \ldots e_n$ are distinct.  Thus $F$ satisfies Condition~(K).  
\end{proof}

\begin{lem}
\label{lem-cofinalplus}

Let $E$ be a graph and let $F$ be a desingularization of $E$.  Then the
following are equivalent:
\begin{itemize}
\item [(1)] $F$ is cofinal;
\item [(2)] $E$ is cofinal and for every singular vertex $v_0 \in E^0$ we have
$E^0 \geq v_0$.
\end{itemize}
\end{lem}

\begin{proof}
Assume $F$ is cofinal and fix $v \in E^0$.  Suppose $\lambda \in
E^\infty$.  Because $F$ is cofinal, $v \geq
\phi_\infty(\lambda)$ in $F$. Thus by Lemma \ref{lem-correspondence}(b), $v
\geq \lambda$ in 
$E$.  Now let $v_0 \in E^0$ be any singular vertex.  Then $\phi_\infty(v_0)$
is the infinite tail $e_1e_2 \ldots$ added to $v_0$.  By cofinality of $F$,
$v$ connects to $e_1e_2 \ldots$, and since any path that connects to
$e_1e_2 \ldots$ connects to $v_0$, we know that there is a path $\alpha \in
F^*$ from $v$ to $v_0$.  But then $\phi^{-1}(\alpha)$ is a path from $v$
to $v_0$ in $E$.  Hence $E^0 \geq v_0$.  

Now assume $E$ is cofinal and for every singular vertex $v_0$ we have $E^0 \geq
v_0$.  If $E$ has a sink $v_0$, then since $E$ is cofinal it follows
that $E^\infty = \emptyset$.  Furthermore, since $E^0 \geq v_0$ it must
be the case that $v_0$ is the only sink in $E$.  Hence $F$ is obtained
from $E$ by adding a single tail at $v_0$.  Now if $\lambda \in
F^\infty$, then since $E^\infty = \emptyset$ we must have that
$\lambda$ eventually ends up in the tail.   If $w \in F^0$, then either
$w$ is in the tail or $w \in E^0$.  Since $E^0 \geq v_0$ this implies
that in either case $w \geq \lambda$.  Hence $F$ is cofinal. 

Now assume that $E$ has no sinks.
Let $\lambda \in F^\infty$ and $v \in F^0$.  We must show that $v \geq
\lambda$ in $F$.  We will first show that it suffices to prove this for the
case when $s(\lambda) \in E^0$ and $v \in E^0$.  If
$v = v_n$, a vertex in one of the added infinite tails, then because $E$ has no
sinks, $v_n$ must be the source of some edge $f_j$ with $r(f_j) \in E^0$ and
 we see that $r(f_j) \geq \lambda$ in $F$ implies $v_n \geq \lambda$ in $F$. 
Likewise, if $s(\lambda) = v_n$, a vertex in the infinite tail added to $v_0$,
then $v \geq e_1 e_2 \ldots e_n \lambda$ in $F$ implies $v \geq \lambda$ in
$F$.  Thus we may replace $\lambda$ by $e_1 e_2 \ldots e_n \lambda$.  Hence
we may assume that $s(\lambda) \in E^0$ and $v \in E^0$. 

Since $\lambda$ is a finite path in $F$ whose source is in $E^0$,
Lemma~\ref{lem-correspondence}(a) implies that $\lambda = \phi_\infty(\mu)$,
where $\mu$ is either an infinite path in $E$ or a finite path in $E$ ending at
a singular vertex.  If $\mu$ is an infinite path, then cofinality of $E$
implies that $v \geq \mu$ and Lemma~\ref{lem-correspondence}(b) implies that $v
\geq \phi_\infty(\mu) = \lambda$.  If $\mu$ is a finite path ending at a
singular vertex, then $v \geq \mu$ by assumption and so $v \geq \phi_\infty(\mu)
= \lambda$.  Thus $F$ is cofinal.
\end{proof}

The next two lemmas will be used to prove Theorem~\ref{thm-morita}, which
states that $C^*(E)$ is isomorphic to a full corner of $C^*(F)$. 
Lemma~\ref{lem-CKEinF} says, roughly speaking, that a Cuntz-Krieger
$F$-family contains a Cuntz-Krieger $E$-family; and
Lemma~\ref{lem-extendEfamtoF} says that we can extend a Cuntz-Krieger
$E$-family to obtain a Cuntz-Krieger $F$-family.  

\begin{lem}
\label{lem-CKEinF}

Suppose $E$ is a graph and let $F$ be a
desingularization of $E$.  
If $\{T_e, Q_v\}$ is a Cuntz-Krieger $F$-family, then there exists a 
Cuntz-Krieger $E$-family in $C^*(\{T_e, Q_v\})$.  

\end{lem}

\begin{proof} 
For every vertex $v$ in $E$, define $P_v := Q_v$.  For every edge $e$ in
$E$ with $s(e)$ not a singular vertex, define $S_e := T_e$.  If $e$ is
an edge in $E$ with $s(e) = v_0$ a singular vertex, then $e = g_j$ for some
$j$, and we define $S_e := T_{\alpha^j}$.   The fact that $\{S_e, P_v\,|\,e
\in E^1, v \in E^0\}$ is a Cuntz-Krieger $E$-family follows immediately from
the fact that $\{T_e, Q_v \,|\, e \in F^1, v \in F^0 \}$ is a Cuntz-Krieger
$F$-family.
\end{proof}

\begin{lem}
\label{lem-extendEfamtoF}

Let $E$ be a graph and let $F$ be 
a desingularization of $E$. For every
Cuntz-Krieger $E$-family $\{S_e, P_v \,|\, e \in E^1, v \in E^0\}$ on
a Hilbert space $\mathcal{H}_E$, there exists a Hilbert space 
$\mathcal{H}_F = \mathcal{H}_E \oplus \mathcal{H}_T$ and a  
Cuntz-Krieger $F$-family $\{T_e, Q_v \,|\, e \in F^1, 
v \in F^0\}$ on $\mathcal{H}_F$ satisfying: 
\begin{itemize}
\item $P_v = Q_v$ for every $v \in E^0$;
\item $S_e = T_e$ for every $e \in E^1$ such that $s(e)$ is not a
singular vertex;
\item $S_e = T_{\alpha^j}$ for every $e = g_j \in E^1$
such that $s(e)$ is a singular vertex;
\item $\sum_{v \notin E^0} Q_{v}$ is the projection onto
$\mathcal{H}_T$.
\end{itemize}
\end{lem}

\begin{proof}
We prove the case where $E$ has just one singular vertex $v_0$.  If $v_0$ is a
sink, then the result follows from \cite[Lemma~1.2]{bprs}.   Therefore let
us assume that $v_0$ is an infinite-emitter.  Given a Cuntz-Krieger $E$-family
$\{S_e, P_v\}$ we define $R_0 := 0$ and $R_n := \sum_{j=1}^n
S_{g_j}S_{g_j}^*$ for each positive integer $n$.  Note that the $R_n$'s are
projections because the $S_{g_j}$'s have orthogonal ranges.  Furthermore, $R_n
\leq R_{n+1} < P_{v_0}$ for every $n$.

Now for every integer $n \geq 1$ define $\mathcal{H}_n
:= (P_{v_0} - R_n)\mathcal{H}_E$ and set  
$$
\mathcal{H}_F := \mathcal{H}_E \oplus \bigoplus_{n=1}^{\infty}
\mathcal{H}_n.
$$

For every $v \in E^0$ define $Q_v = P_v$ acting on the $\mathcal{H}_E$
component of $\mathcal{H}_F$ and zero elsewhere.  That is, $Q_v(\xi_E, \xi_1,
\xi_2, \dots) = (P_v\xi_E, 0, 0, \dots)$.  Similarly, for every $e \in E^1$
with $s(e) \neq v_0$ define $T_e = S_e$ on the $\mathcal{H}_E$ component.
For each vertex $v_n$ on the infinite tail define $Q_{v_n}$ to
be the projection onto $\mathcal{H}_n$.  That is, $Q_{v_n}(\xi_E,
\xi_1, \ldots ,\xi_n, \ldots) = (0,0,\ldots,\xi_n,0,\ldots)$.  
Now note that because the $R_n$'s are non-decreasing,
$\mathcal{H}_{n} \subseteq
\mathcal{H}_{n-1}$ for each $n$.  Thus
for each edge of the form $e_n$
we can define $T_{e_n}$ to be the inclusion of $\mathcal{H}_n$ into
$\mathcal{H}_{n-1}$ (where $\mathcal{H}_0$ is taken to mean
$P_{v_0}\mathcal{H}_E$).  More precisely, 
$$T_{e_n}(\xi_E, \xi_1, \xi_2, \dots) =
(0,0,\dots,0,\xi_n,0,\dots),
$$ where the $\xi_n$ is in the $\mathcal{H}_{n-1}$ component.  

Finally, for each edge $g_j$ and for each $\xi \in \mathcal{H}_E$ we
have $S_{g_j}\xi \in \mathcal{H}_{j-1}$.  To see this recall
that $\mathcal{H}_{j-1} = (P_{v_0} - R_{j-1})\mathcal{H}_E$, and thus
$(P_{v_0} - R_{j-1})S_{g_j}\xi = S_{g_j}\xi$.  Therefore we can define
$T_{f_j}$ by
$$
T_{f_j}(\xi_E, \xi_1, \xi_2, \dots) = (0,\dots,0,S_{g_j}\xi_E,0,\dots),
$$
where the nonzero term appears in the $\mathcal{H}_{j-1}$ component.

We will now check that the collection $\{T_e, Q_v\}$ is a Cuntz-Krieger
$F$-family.  It follows immediately from definitions and the
Cuntz-Krieger relations on $E$ that $T_e^*T_e = Q_{r(e)}$ for every
$e$ that is not of the form $f_j$ or $e_n$, and that $Q_v =
\sum_{s(e)=v} T_eT_e^*$ for every $v$ not on the infinite tail.  Furthermore, it
is easy to check using the definitions that the $Q_v$'s are mutually orthogonal
and that $T_{e_n}^*T_{e_n} = Q_{r(e_n)}$ for every edge $e_n$ on the
infinite tail.  Now note that for every $f_j$, 
\begin{align*}
T_{f_j}^*T_{f_j}(\xi_E, \xi_1, \xi_2, \dots) &=
T_{f_j}^*(0,\dots,0,S_{g_j}\xi_E,0,\dots)\\ 
&= (S_{g_j}^*S_{g_j}\xi_E, 0, 0, \dots)\\ 
&= (P_{r(e_j)}\xi_E, 0, 0, \dots)\\
&= Q_{r(e_j)}(\xi_E, \xi_1, \xi_2, \dots).
\end{align*}

Finally, let $v_n$ be a vertex on the infinite tail.  The edges
emanating from $v_n$ are $e_{n+1}$ and $f_{n+1}$, and we have 
$$
T_{e_{n+1}}T_{e_{n+1}}^*(\xi_E, \xi_1, \dots) =
(0,\dots,0,(P_{v_0}-R_{n+1})\xi_n,0,\dots),
$$
where the nonzero term is in the $\mathcal{H}_n$ component.  Also 
$$
T_{f_{n+1}}T_{f_{n+1}}^*(\xi_E, \xi_1, \dots) =
(0,\dots,0,S_{g_{n+1}}S_{g_{n+1}}^*\xi_n,0,\dots),
$$
where the nonzero term is again in the $\mathcal{H}_n$ component.  We then have
the following:
\begin{align*}
(T_{e_{n+1}}T_{e_{n+1}}^* + T_{f_{n+1}}T_{f_{n+1}}^*) (\xi_E, \xi_1,
\dots) 
&= (0,\dots,0,(P_{v_0} - R_{n+1} + S_{g_{n+1}}S_{g_{n+1}}^*)\xi_n, 0, \dots)\\
&= (0,\dots,0,(P_{v_0} - R_n)\xi_n, 0, \dots)\\
&= (0,\dots,0,\xi_n,0,\dots)\\
&= Q_{v_n}(\xi_E, \xi_1, \dots).
\end{align*}
Thus $\sum_{ \{e : s(e) = v_n \}} T_eT_e^* = T_{e_{n+1}}T_{e_{n+1}}^* +
T_{f_{n+1}}T_{f_{n+1}}^* = Q_{v_n} = Q_{r(e_n)}$ and we have established that
$\{T_e, Q_v\}$ is a Cuntz-Krieger $F$-family.  It is easy to verify that the
bulleted points in the statement of the lemma are satisfied.  
\end{proof}

\begin{thm}

\label{thm-morita}
Let $E$ be a graph and let $F$ be a
desingularization of $E$.  Then $C^*(E)$ is isomorphic to a full
corner of $C^*(F)$.  Consequently, $C^*(E)$ and $C^*(F)$ are Morita
equivalent.

\end{thm}

\begin{proof}
Again for simplicity we assume that $E$ has only one singular vertex $v_0$. Let
$\{t_e, q_v \,|\, e \in F^1, v \in F^0\}$ denote the canonical set of
generators for $C^*(F)$ and let $\{s_e, p_v \,|\, e \in E^1, v \in
E^0\}$ denote the Cuntz-Krieger $E$-family in $C^*(F)$ constructed in
Lemma~\ref{lem-CKEinF}.  Define $B:=C^*(\{s_e, p_v\})$ and 
$p:=\sum_{v \in E^0} q_v$.  To prove
the proposition, we will show that $C^*(E) \cong B \cong pC^*(F)p$ is a full
corner in $C^*(F)$.

Since $B$ is generated by a Cuntz-Krieger $E$-family, in order to
show that $B \cong C^*(E)$ it suffices to prove that $B$ satisfies
the universal property of $C^*(E)$. 
Let $\{S_e, P_v \,|\, e \in E^1,
v \in E^0\}$ be a Cuntz-Krieger $E$-family on a Hilbert space
$\mathcal{H}_E$.  Then by Lemma \ref{lem-extendEfamtoF} we can
construct a Hilbert space $\mathcal{H}_F$ and a Cuntz-Krieger
$F$-family $\{T_e, Q_v \,|\, e \in F^1, v \in F^0\}$ on
$\mathcal{H}_F$ such that $Q_v = P_v$ for every $v \in E^0$, $T_e =
S_e$ for every $e \in F^1$ with $s(e) \neq v_0$, and $S_{g_j} =
T_{\alpha^j}$ for every edge $g_j$ in $E$ whose source is $v_0$.
Now by the universal property of $C^*(F)$, we have a homomorphism
$\pi$ from $C^*(F)$ onto $C^*(\{T_e, Q_v \,|\, e \in F^1, v \in
F^0\})$ that takes $t_e$ to $T_e$ and $q_v$ to $Q_v$.   

For any $v \in E^0$ we have $p_v = q_v$, so $\pi(p_v) = Q_v = P_v$.  Let $e
\in E^1$.  If $s(e) \neq v_0$, then $s_e = t_e$ and $\pi(s_e) = T_e =
S_e$.  Finally, if $s(e) = v_0$ then $e = g_j$ for some $j$, and $s_e =
t_{\alpha^j}$ so that $\pi(s_{g_j}) = T_{\alpha^j} =
S_{g_j}$.  Thus $\pi|_B$ is a representation of $B$ on $\mathcal{H}_E$ that
takes generators of $B$ to the corresponding elements of
the given Cuntz-Krieger $E$-family.  Therefore $B$ satisfies the
universal property of $C^*(E)$ and $C^*(E) \cong B$.

We now show that $B \cong pC^*(F)p$.  Just as in \cite[Lemma 1.2(c)]{bprs},
we have that $\sum_{v \in E^0} q_v$ converges strictly in $M(C^*(F))$ to a
projection $p$ and that for any $\mu, \nu \in F^*$ with $r(\mu)=r(\nu)$,  
\begin{equation}
\label{eq-pp}
p t_{\mu} t_{\nu}^* p = \left\{\begin{array}{ll}
t_{\mu} t_{\nu}^* & \hbox{if $s(\mu), s(\nu) \in E^0$;}\\
0 & \hbox{otherwise.}\\
\end{array}\right.
\end{equation}
Therefore the generators of $B$ are contained in $pC^*(F)p$ and $B
\subseteq pC^*(F)p$. To show the reverse inclusion, let $\mu$ and $\nu$ be
finite paths in $F$ with $r(\mu) = r(\nu)$.  We need to show that 
$pt_{\mu} t_{\nu}^*p 
\in B$.  If either $\mu$ or $\nu$ does not start in
$E^0$, then $pt_{\mu} t_{\nu}^*p = 0$ by the above formula.  Hence we 
may as well assume that both $\mu$ and $\nu$ start in $E^0$.  Now, if $r(\mu) =
r(\nu) \in E^0$ as well, there will exist unique $\mu', \nu' \in E^0$ with
$\phi(\mu') = \mu$ and $\phi(\nu') = \nu$.  In this case, $t_\mu =
s_{\mu'}$ and $t_\nu = s_{\nu'}$, so 
$$
pt_{\mu} t_{\nu}^*p = t_{\mu}t_{\nu}^* = s_{\mu'}s_{\nu'}^* \in B. 
$$
On the other hand, if $r(\mu) = r(\nu) \not\in E^0$ then $r(\mu) = r(\nu) =
v_n$ for some $n$.  We shall prove that $pt_{\mu} t_{\nu}^*p \in B$ by
induction on $n$.  Suppose that $pt_{\mu'} t_{\nu'}^*p \in B$ for any paths
$\mu'$ and $\nu'$ with $r(\mu')=r(\nu')=v_{n-1}$.  Then if
$r(\mu)=r(\nu)=v_n$ we shall write $\mu = \mu' e_n, \nu =\nu' e_n$ for
finite paths $\mu'$ and $\nu'$ with $r(\mu')=r(\nu')=v_{n-1}$.   Now there are
precisely two edges, $e_n$ and $f_{n}$ with source $v_{n-1}$.  Thus
\begin{align*}
pt_{\mu}t_{\nu}^*p 
&= pt_{\mu'}t_{e_n}t_{e_n}^* t_{\nu'}^*p\\
&= pt_{\mu'} (q_{v_{n-1}} - t_{f_{n}}t_{f_{n}}^*) t_{\nu'}^*p\\
&= pt_{\mu'} t_{\nu'}^*p - pt_{\mu'f_{n}}t_{\nu'f_{n}}^*p,
\end{align*}
which is in $B$.  Hence $pC^*(F)p \subseteq B$.  Finally, we
note that $pC^*(F)p$ is full by an argument identical to the one given
in \cite[Lemma 1.2(c)]{bprs}.
\end{proof}

Theorem \ref{thm-morita} allows us to get easy proofs of several known results
by passing to a desingularization and using the corresponding result for
row-finite graphs.

\begin{cor}
\label{cor-uniqueness}

Suppose $E$ is a graph in which every loop has an exit, and that $\{S_e, 
P_v\}$ and $\{T_e, Q_v\}$ are two Cuntz-Krieger $E$-families in which all
the projections $P_v$ and $Q_v$ are non-zero.  Then there is an isomorphism
$\phi:C^*(\{S_e, P_v\}) \rightarrow C^*(\{T_e, Q_v\})$ such that $\phi(S_e)
= T_e$ for all $e \in E^1$
and $\phi(P_v) = Q_v$ for all $v \in E^0$.

\end{cor}

\begin{proof}
Let $F$ be a desingularization of $E$.
Use Lemma \ref{lem-extendEfamtoF} to construct $F$-families from
the given $E$-families.  Then apply \cite[Theorem 3.1]{bprs} to get an
isomorphism between the $C^*$-algebras generated by the $F$-families that will
restrict to an isomorphism between $C^*(\{S_e, P_v\})$ and $C^*(\{T_e,Q_v\})$.
\end{proof}

\begin{cor}
\label{cor-AF}

Let $E$ be a graph.  Then $C^*(E)$ is an AF-algebra if and only if $E$ has
no loops.  

\end{cor}

\begin{proof}
This follows from \cite[Theorem 2.4]{kpr}
and the fact that the class of AF-algebras is closed under stable
isomorphism (see \cite[Theorem 9.4]{effros}).
\end{proof}

\begin{cor}
\label{cor-pureinf}

Let $E$ be a graph. Then $C^*(E)$ is purely infinite if and only if every
vertex in $E$ connects to a loop and every loop in $E$ has an exit.

\end{cor}

\begin{proof}
By \cite[Proposition 5.3]{bprs} and the fact that pure infiniteness is
preserved by passing to corners, every vertex connects to a loop and every
loop has an exit implies pure infiniteness.  For the converse we note that
the proof given in \cite[Theorem 3.9]{kpr} works for arbitrary graphs.
\end{proof}

The following result generalizes 
\cite[Theorem 3]{flr} 
and \cite[Corollary 4.5]{fr} and it was proven independently in \cite{sz} and
\cite{patgg}.

\begin{cor}
\label{cor-simple}

Let $E$ be a graph.  Then $C^*(E)$ is simple if and only if
\begin{itemize}
\item [(1)] every loop in $E$ has an exit;
\item [(2)] $E$ is cofinal;
\item [(3)] for every singular vertex $v_0 \in E^0$, $E^0 \geq v_0$.
\end{itemize}

\end{cor}

\begin{proof}
Letting $F$ denote a desingularization of $E$, we have
$$
\begin{array}{lcll}
C^*(E) \hbox{ is simple} & \Longleftrightarrow & C^*(F) \hbox{ is
simple} & \hbox{(by Theorem \ref{thm-morita})}\\
& \Longleftrightarrow & F \hbox{ is cofinal and every loop in $F$ has an
exit} &
\hbox{(by \cite[Proposition 5.1]{bprs})}\\
& \Longleftrightarrow & E \hbox{ satisfies (1),(2), and (3)} &
\hbox{(by Lemmas \ref{lem-LK} and \ref{lem-cofinalplus})}.\\
\end{array}
$$
\end{proof}

\begin{rem}
\label{rem-dichotomy}
We see from the results above that the dichotomy described in
\cite[Remark 5.6]{bprs} holds for arbitrary graphs:  If $C^*(E)$ is simple, then
it is either AF or purely infinite.  For if $E$ has no loops then
Corollary~\ref{cor-AF} shows that $C^*(E)$ is AF.  If $E$ does have loops, then
Corollary~\ref{cor-simple} says that they all have exits and that $E$ is
cofinal;  thus every vertex connects to every loop and
Corollary~\ref{cor-pureinf} applies.
\end{rem}

%------------------------------------------------------------------------
\section{Ideal structure}
\label{sec-ideal}
%-----------------------------------------------------------------------

Let $E$ be a directed graph.
A set $H \subseteq E^0$ is {\it hereditary} if whenever $v \in H$ and
$v \geq w$, then $w \in H$.  A hereditary set $H$ is called {\it
saturated} if every vertex that is not a singular vertex and that
feeds only into $H$ is itself in $H$;  that is, if 
$$
\hbox{$v$ not singular and $\{r(e) \,|\, s(e)=v\} \subseteq H$ 
implies $v \in H$.}
$$
If $E$ is row-finite this definition reduces to the one given in
\cite{bprs}.  It was shown in \cite[Theorem 4.4]{bprs} that 
if $E$ is row-finite and satisfies Condition~(K), then
every saturated hereditary subset $H$
of $E^0$ gives rise to exactly one ideal $\hbox{$I_H :=$ the ideal
generated by $\{p_v \,|\, v \in H\}$}$ in $C^*(E)$.  
If $E$ is a graph that is not row-finite, it is easy to check that with
the above definition of saturated \cite[Lemma 4.2]{bprs} and
\cite[Lemma 4.3]{bprs} still hold.  Consequently, $H
\mapsto I_H$ is still injective, just as in the proof of \cite[Theorem
4.1]{bprs}.  However, it is no longer true that this map is surjective;
that is, there may exist ideals in $C^*(E)$ that are not of the form $I_H$
for some saturated hereditary set $H$.  The reason the proof for row-finite
graphs no longer works is that if $I$ is an ideal, then $\{s_e + I, p_v
+ I\}$ will not necessarily be a Cuntz-Krieger $E \setminus H$-family for
the graph $E \setminus H$ defined in \cite[Theorem~4.1]{bprs}.  
It turns out that to describe an
arbitrary ideal in $C^*(E)$ we need a saturated hereditary subset and
one other descriptor.  Loosely speaking, this descriptor tells us how close
$\{s_e + I, p_v + I\}$ is to being a Cuntz-Krieger $E \setminus H$-family.

Given a saturated hereditary subset $H
\subseteq E^0$, define 
$$
B_H := \{v \in E^0 \,|\, \hbox{$v$ is an infinite-emitter and $0 <
|s^{-1}(v) \cap r^{-1}(E^0 \setminus H)| < \infty$}\}.
$$
Therefore $B_H$ is the set of infinite-emitters that point to only a finite
number of vertices not in $H$.   Since $H$ is hereditary, $B_H$ is disjoint
from $H$.   Now fix a saturated hereditary subset $H$ and let $S$ be any
subset of $B_H$.  Let $\{s_e, p_v\}$ be the canonical generating Cuntz-Krieger
$E$-family and define
$$
I_{(H,S)} := \hbox{the ideal in $C^*(E)$ generated by $\{p_v \,|\, v
\in
H\} \cup
\{\pvoh \,|\, v_0 \in S\}$},
$$
where 
$$
\pvoh := p_{v_0} - \sum_{{s(e) = v_0} \atop {r(e) \notin H}} s_e s_e^*. 
$$
Note that the definition of $B_H$ ensures that the sum on the right is
finite.

Our goal is to show  that the correspondence $(H,S) \mapsto I_{(H,S)}$ is
a lattice isomorphism, so we must describe the lattice structure on 
$$
\{(H,S) \,|\,
\hbox{$H$ is a saturated hereditary subset of $E^0$ and $S \subseteq 
B_H$}\}.
$$  
We say $(H,S) \leq (H',S')$ if and only if $H \subseteq H'$
and $S \subseteq H' \cup S'$.  With this definition, the reader who
is willing to spend a few minutes can
check using nothing more than basic set theory that the
following equations define a greatest lower bound and
least upper bound:
$$
(H_1,S_1) \wedge (H_2,S_2) := 
\left((H_1 \cap H_2), \left(H_1 \cup H_2 \cup S_1 \cup S_2)\right)
\cap B_{H_1 \cap H_2} \right)
$$
$$
(H_1,S_1) \vee (H_2,S_2) := 
\left(
\bigcup_{n=0}^\infty X_n \ , \
(S_1 \cup S_2) \cap B_{ \bigcup_{n=0}^\infty X_n } 
\right)$$
where $X_n$ is defined recursively as $X_0 := H_1 \cup H_2$ and $X_{n+1} :=
X_n \cup \{ v \in E^0 \,|\, 0 < | s^{-1}(v) | < \infty \text{ and } \{ r(e)
\,|\, s(e) = v \} \subseteq X_n \} \cup \{ v \in E^0 \,|\, v \in S_1 \cup S_2
\text{ and } \{r(e) \,|\, s(e)=v \} \subseteq X_n \}$.  The reason for this
strange definition of the $X_n$'s is the following:  If $Y_0$ is a hereditary
subset, then the saturation of $Y_0$ may be defined as the increasing union
of $Y_{n+1} := Y_n \cup \{ v \in E^0 \,|\, 0 < | s^{-1}(v) | < \infty
\text{ and } \{ r(e) \,|\, s(e) = v \} \subseteq Y_n \}$.  In the $X_n$'s
above we need not only these elements, but also at each stage we must include
the infinite emitters in $S_1 \cup S_2$ that only feed into $X_n$.

We now describe a correspondence between pairs $(H,S)$ as above and saturated
hereditary subsets of vertices in a desingularization of $E$.  
Suppose that $E$ is a graph and let $F$ be a desingularization of $E$.
Also let $H$ be a
saturated hereditary subset of $E^0$ and let $S \subseteq B_H$.  We
define a saturated hereditary subset $H_S \subseteq F^0$.  First set
$\Htil := H \cup \{v_n \in F^0 \,|\, v_n \hbox{ is on a tail added
to a vertex in $H$}\}$.  Now for each $v_0 \in S$ let $N_{v_0}$ be the
smallest nonnegative integer such that $r(e_j) \in
H$ for all $j \geq N_{v_0}$.  The number $N_{v_0}$ exists since
$v_0 \in B_H$ implies that there must be a vertex on the tail added to $v_0$
beyond which each vertex points only to the next vertex on the tail and 
into $H$. Define
$T_{v_0} := \{v_n \,|\, v_n \hbox{ is on the infinite tail added to 
$v_0$ and $n \geq N_{v_0}$}\}$ and define
$$
H_S := \Htil \cup \bigcup_{v_0 \in S} T_{v_0}. 
$$
Note that for $v_0 \in B_H$ we have $v_0 \not \in H_S$.  Furthermore, the
tail attached to $v_0$ will eventually be inside $H_S$ if and only if $v_0
\in S$. It is easy to check that $H_S$ is hereditary, and choosing $N_{v_0}$
to be minimal ensures that $H_S$ is saturated.

\begin{ex}
\label{ex-prim}

Suppose $E$ is the following graph:
$$
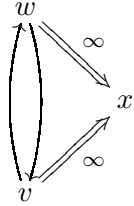
\xymatrix{
w \ar@/^/[dd] \ar@{=>}[dr]^{\infty} & \\
& x \\
v \ar@/^/[uu] \ar@{=>}[ur]_{\infty}\\
}
$$
A desingularization $F$ is given by
$$
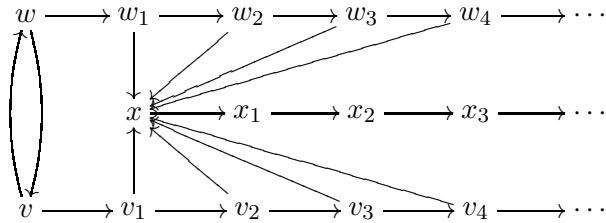
\xymatrix{
w \ar@/^/[dd] \ar[r] & w_1 \ar[r] \ar[d] & w_2
\ar[r] \ar[dl] & w_3 \ar[r] \ar[dll] & w_4 \ar[r] \ar[dlll] & \cdots\\
& x \ar[r] & x_1 \ar[r] & x_2 \ar[r] & x_3 \ar[r] & \cdots\\
v \ar@/^/[uu] \ar[r] & v_1 \ar[r] \ar[u] & v_2
\ar[r] \ar[ul] & v_3 \ar[r] \ar[ull] & v_4 \ar[r] \ar[ulll] & \cdots\\
}
$$

The only saturated hereditary (proper) subset in $E$ is the set $H 
= \{x\}$.  In this case $B_H = \{v,w\}$.  There are four subsets of
$B_H$ and there are four
saturated hereditary (proper) subsets in the desingularization.  In
particular, if
$S = \emptyset$, then $H_S$ consists of only the tail added to $x$;  if $S$
contains $w$, then $H_S$ also includes $\{w_2, w_3, \dots\}$; and if $S$
contains $v$, then $H_S$ also includes $\{v_2, v_3, \dots\}$. 

\end{ex}

The proof of the following lemma is straightforward.

\begin{lem}
\label{lem-HSHS}
Let $E$ be a graph and let $F$ be a desingularization of $E$.
The map $(H,S) \mapsto H_S$ is an isomorphism from the lattice 
$$
\{(H,S) \,|\,
\hbox{$H$ is a saturated hereditary subset of $E^0$ and $S \subseteq 
B_H$}\}
$$  
onto the lattice of saturated hereditary subsets of $F$.
\end{lem}

Suppose $E$ is a graph that satisfies Condition~(K) and $F$ is a
desingularization of $E$. 
Because $C^*(E)$ is isomorphic to the full corner
$pC^*(F)p$, we have that $C^*(E)$ and $C^*(F)$ are Morita equivalent via the
imprimitivity bimodule $pC^*(F)$.   It then follows from \cite[Proposition
3.24]{tfb} that the Rieffel correspondence between ideals in $C^*(F)$ and
ideals in $C^*(E)$ is given by the map $I \mapsto pIp$.  

\begin{prop}
\label{prop-HSHS}

Let $E$ be a graph satisfying Condition~(K) and let $F$ be a
desingularization of $E$. Let $H$ be a
saturated hereditary subset of $E^0$ and let $S \subseteq B_H$.  If
$\{t_e, q_v\}$ is a generating Cuntz-Krieger $F$-family and $p 
= \sum_{v \in E^0} q_v$, then $pI_{H_S}p = I_{(H,S)}$.

\end{prop}

\begin{proof}
That $pI_{H_S}p \subseteq I_{(H,S)}$ is immediate from (\ref{eq-pp}).
We show the reverse inclusion by showing that the generators of
$I_{(H,S)}$ are in $pI_{H_S}p$.  Letting $\{s_e, p_v\}$ denote the
Cuntz-Krieger $E$-family defined in the proof of Lemma~\ref{lem-CKEinF}, the
generators for $I_{(H,S)}$ are $\{p_v \,|\, v \in H\} \cup \{p_{v_0}^H \,|\,
v_0 \in S\}$. Clearly for $v \in H$, we have $p_v = q_v
= pq_vp \in pI_{H_S}p$, so all that remains to show is that for every $v_0 \in
S$ we have $\pvoh \in pI_{H_S}p$.  

Let $v_0 \in S$ and $n:= N_{v_0}$.  Then 
\begin{align*}
q_{v_0} &= t_{e_1}t_{e_1}^* + t_{f_1}t_{f_1}^* \\
&= t_{e_1}q_{v_1}t_{e_1}^* + t_{f_1}t_{f_1}^* \\
&= t_{e_1}(t_{e_2}t_{e_2}^*+t_{f_2}t_{f_2}^*)t_{e_1}^* + t_{f_1}t_{f_1}^*\\
&= t_{e_1e_2}q_{v_2}t_{e_1e_2}^* + t_{e_1f_2}t_{e_1f_2}^* +t_{f_1}t_{f_1}^*\\
& \quad \quad \quad \vdots \\
&= t_{e_1 \ldots e_n}t_{e_1 \ldots e_n}^* + \sum_{j=1}^n
t_{\alpha^j}t_{\alpha^j}^* 
\end{align*}
Now since $r(e_n) = v_n \in H_S$ we see that $q_{v_n} \in I_{H_S}$ and hence
$t_{e_n} = t_{e_n}t_{e_n}^*t_{e_n} = t_{e_n}q_{v_n} \in I_{H_S}$. 
Consequently, $t_{e_1 \ldots e_n}t_{e_1 \ldots e_n}^* \in I_{H_S}$. 
Similarly, whenever $r(\alpha^j) \in H$, then $t_{\alpha^j}t_{\alpha^j}^* \in
I_{H_S}$.  Now, by definition, every $\alpha^j$ with $r(\alpha^j) \notin H$
has $j < n$.  Therefore the above equation shows us that
\begin{align*}
\pvoh &= p_{v_0} - \sum_{{s(g_j) = v_0} \atop {r(g_j) \notin H}} s_{g_j}
s_{g_j}^* \\
&= q_{v_0} - \sum_{{s(\alpha^j) = v_0} \atop {r(\alpha^j) \notin H}}
t_{\alpha^j} t_{\alpha^j}^* \\
&= \sum_{ {r(\alpha^j) \in H} \atop {j < n} }
t_{\alpha^j} t_{\alpha^j}^* + t_{e_1 \ldots e_n} t_{e_1 \ldots e_n}^*
\end{align*}
which is an element of $I_{H_S}$ by the previous paragraph.  Hence $I_{H_S}
\subseteq I_{H,S}$.
\end{proof}

\begin{cor}
\label{cor-prim}
Let $E$ be a graph satisfying Condition~(K) and let $F$ be a
desingularization of $E$. If $H$ is a
saturated hereditary subset of $E^0$ and $S \subseteq B_H$, then $I_{(H,S)}$
is a primitive ideal in $C^*(E)$ if and only if $I_{H_S}$ is a primitive
ideal in $C^*(F)$.
\end{cor}

We now have the following:
$$
\xymatrix{
\{(H,S) \,|\, \hbox{$H$ is saturated, hereditary in $E$ and $S \subseteq
B_H$}\} \ar[d] \ar@{.>}[rr] &&
\hbox{ideals in $C^*(E)$}\\
\hbox{saturated, hereditary subsets of $F$} \ar[rr] && \hbox{ideals in
$C^*(F)$.} \ar[u] 
}
$$
The map on the left is $(H,S) \mapsto H_S$, which is a lattice isomorphism by
Lemma \ref{lem-HSHS}.  The lattice isomorphism $H \mapsto I_H$ across the
bottom comes from \cite[Theorem~4.4]{bprs}.  The map on the right is $I_{H_S}
\mapsto I_{(H,S)}$ and is an isomorphism because it agrees with the Rieffel
correspondence (Proposition \ref{prop-HSHS}).  Composing the three yields the
following:

\begin{thm}
\label{thm-ideal}

Let $E$ be a graph that satisfies Condition~(K).  Then the map
$(H,S) \mapsto I_{(H,S)}$ is a lattice isomorphism from the lattice 
$$
\{(H,S) \,|\,
\hbox{$H$ is a saturated hereditary subset of $E^0$ and $S \subseteq
B_H$}\}
$$ 
onto the lattice of ideals in $C^*(E)$.

\end{thm}

%------------------------------------------------------------------------
\section{Primitive ideal space}
\label{sec-prim}
%-----------------------------------------------------------------------

The following definition generalizes that in \cite[Proposition~6.1]{bprs}.

\begin{defn}
\label{defn-maxtail}
Let $E$ be a graph.  A nonempty subset $\gamma \subseteq E^0$ 
is called a {\it maximal tail} if it satisfies the following conditions:
\begin{itemize}
\item [(a)] for every $w_1, w_2 \in \gamma$ there exists $z \in \gamma$
such that $w_1 \geq z$ and $w_2 \geq z$;  
\item [(b)] for every $v \in \gamma$ that is not a singular vertex, there
exists an edge $e$ with $s(e)=v$ and $r(e) \in \gamma$;
\item [(c)] $v \geq w$ and $w \in \gamma$ imply $v \in \gamma$.
\end{itemize}
\end{defn}

Given a graph $E$ we denote by $\Lambda_E$ the set of all maximal tails in
$E$.  Note that if $v_0$ is a sink, then the set $\lambda_{v_0} := \{v \in E^0
\,|\, v \geq v_0\}$ is a maximal tail according to
Definition~\ref{defn-maxtail}, but was not considered to be a maximal tail
in \cite[Section~6]{bprs}.  In addition, when $v_0$ is an infinite-emitter
$\lambda_{v_0} := \{v \in E^0 \,|\, v \geq v_0\}$ is a maximal
tail. 

\begin{defn}
\label{break-vertex}
If $E$ is a graph, then a \emph{breaking vertex} is an element $v \in E^0$
such that $| s^{-1}(v) | = \infty$ and $0 < | \{ e \in E^1 \, | \, s(e)=v
\text{ and } r(e) \geq v \} | < \infty$.  We denote the set of breaking
vertices of $E$ by $BV(E)$.
\end{defn}

\begin{rem}
Notice that if $H$ is a hereditary subset in a graph $E$ and $v_0 \in B_H$,
then $v_0$ is a breaking vertex if and only if there exists an edge $e \in
E^1$ with $s(e)=v_0$ and $r(e) \geq v_0$.  Also note that if $H$ is a
saturated hereditary subset in a graph $E$ and $E^0 \setminus H =
\lambda_{v_0}$ for some singular vertex $v_0$, then $v_0 \in B_H$ if and only
if $v_0$ is a breaking vertex.
\end{rem}

We let $\Xi_E := \Lambda_E \cup BV(E)$ denote the disjoint union of the
maximal tails and the breaking vertices.  We shall see that the elements of
$\Xi_E$ correspond to the primitive ideals in $C^*(E)$.

\begin{lem}
\label{tailpath}
If $E$ is a graph and $\gamma$ is a maximal tail in
$E$, then $\gamma = \{ v \in E^0 \, | \, v \geq \alpha \}$ for some $\alpha
\in  E^\infty \cup \{\alpha \in E^* \,|\, r(\alpha) \hbox{ is a
singular vertex}\}$.
\end{lem}

\begin{proof}
It is straightforward to see that if $\alpha \in E^\infty \cup \{\alpha \in
E^* \,|\, r(\alpha) \hbox{ is a singular vertex}\}$, then $\{ v \in E^0
\, | \, v \geq \alpha \}$ is a maximal tail \cite[Remark~6.4]{bprs}. 

Conversely, suppose that $\gamma$ is a maximal tail.  We shall create a path
in $E$ inductively.  Begin with an element $w \in \gamma$.  If there exists
an element $w' \in \gamma$ for which $w' \ngeq w$, then we may use
property~(a) of maximal tails to choose a path $\beta^1$ with $s(\beta^1) =
w$ and $w' \geq r(\beta^1)$.  Now having chosen $\beta^i$, we do one of two
things:  if $w' \geq r(\beta^i)$ for all $w' \in \gamma$, we stop.  If there
exists $w' \in \gamma$ such that $w' \ngeq r(\beta^i)$, then we choose a path
$\beta^{i+1}$ with $s(\beta^{i+1}) = r(\beta^i)$ and $w' \geq
r(\beta^{i+1})$.  We then continue in this manner to produce a path $\beta :=
\beta^1 \beta^2 \ldots$, which may be either finite or infinite.  Note that
since $\gamma$ has either a finite or countable number of elements, we may
choose $\beta$ in such a way that $w \geq \beta$ for all $w \in \gamma$.

Now if $\beta$ is an infinite path we define $\alpha := \beta$.  On the other
hand, if $\beta$ is a finite path then one of two things must occur.  Either
$r(\beta)$ is a singular vertex or there is an edge $e_1 \in E^1$ with $s(e_1)
= r(\beta)$ and $r(e_1) \in \gamma$.  Continuing in this way, we see that
having chosen $e_i$, either $r(e)$ is a singular vertex or there exists
$e_{i+1} \in E^1$ with $s(e_{i+1}) = r(e_i)$ and $r(e_{i+1}) \in \gamma$. 
Using this process we may extend $\beta$ to a path $\alpha:= \beta e_1 e_2
\ldots$ that is either infinite or is finite and ends at a singular vertex.

Now since every vertex on $\alpha$ is an element of $\gamma$ we
certainly have $\{ v \in F^0 \, | \, v \geq \alpha \} \subseteq \gamma$. 
Also, for every element $v \in \gamma$ there exists an $i$ such that $v \geq
r(\beta^i) \geq \alpha$ so we have $\gamma \subseteq \{ v \in F^0 \, | \, v
\geq \alpha \}$.
\end{proof}

\begin{thm}
Let $E$ be a graph.  An ideal $I$ in $C^*(E)$ is a primitive ideal if and only
if one of the following two statements holds:
\begin{enumerate}
\item $I = I_{(H,S)}$, where $E^0 \setminus H$ is a
maximal tail and $S = B_H$; or
\item $I = I_{(H,S)}$, where $E^0 \setminus H =
\lambda_{v_0}$ for some breaking vertex $v_0$ and $S = B_H \setminus \{v_0\}$.
\end{enumerate}
\end{thm}

\begin{proof}
It follows from Theorem~\ref{thm-ideal} that any ideal in $C^*(E)$ has the
form $I_{(H,S)}$ for some saturated hereditary set $H \subseteq E^0$ and some
$S \subseteq B_H$.  Let $F$ be a desingularization of $E$.  It follows from
Corollary~\ref{cor-prim} that $I_{(H,S)}$ is primitive if and only if
$I_{H_S}$ is primitive.

Now suppose that $I_{(H,S)}$, and hence $I_{H_S}$, is primitive.  It follows
from \cite[Proposition~6.1]{bprs} that $F^0 \setminus H_S$ is a maximal tail
in $F$.  Thus by Lemma~\ref{tailpath} we have $F^0 = \{ w \in F^0 \, | \, w
\geq \alpha \}$ for some $\alpha \in F^\infty$.  Now
$\phi_\infty^{-1}(\alpha)$ is either an infinite path in $E$ or a finite path
in $E$ ending at a singular vertex.  In either case $\gamma := \{ w \in E^0
\, | \, w \geq \phi_\infty^{-1}(\alpha) \}$ is a maximal tail in $E$. 
Furthermore, 
$$
v \in E^0 \setminus H 
\Longleftrightarrow
v \notin H
\Longleftrightarrow
v \notin H_S
\Longleftrightarrow
v \geq \alpha \text{ in $F$}
\Longleftrightarrow 
v \geq \phi_{\infty}^{-1}(\alpha) \text{ in $E$}
\Longleftrightarrow 
v \in \gamma.
$$
Therefore $E^0 \setminus H = \gamma$ is a maximal tail.

Now if $S = B_H$, then we are in the case described in part~(1) of the
theorem and the claim holds.  Let us therefore suppose that there exists $v_0
\in B_H \setminus S$.  If we define $T_{v_0} := \{ v_0, v_1, v_2, \ldots \}$
to be the vertices on the tail added to $v_0$, then we see that $v_0 \notin S$
implies that $T_{v_0} \subseteq F^0 \setminus H_S = \{ w \in F^0 \, | \, w
\geq \alpha \}$.  Now for each vertex $v_i$ with $i \geq N_{v_0}$ there are
two edges, $e_{i+1}$ and $f_{i+1}$, with source $v_i$.  Since $r(f_{i+1}) \in
H_S$ and $r(e_{i+1}) = v_{i+1}$, it must be the case that $\alpha$ has the
form $\alpha = \alpha' e_1e_2e_3 \ldots$ for some finite path $\alpha'$ in
$F$.  Consequently, $\phi_\infty^{-1}(\alpha)$ is a finite path in $E$ ending
at $v_0$, and $\gamma = \lambda_{v_0}$.  Now let $X := \{ e \in E^1 \, | \,
s(e)=v_0 \text{ and } r(e) \geq v_0\}$.  Note that if $s(e) = v_0$ and $r(e)
\geq v_0$, then $r(e) \notin H$ since $H$ is hereditary.  Because $v_0 \in
B_H$ it follows that we must have $| X | < \infty$.  Furthermore, since $v_0
\in B_H$ there exists $e \in E^1$ with $s(e)=v_0$ and $r(e) \notin H$.  But
then $r(e) \in \gamma$ and $r(e) \geq \phi_\infty^{-1}(\alpha)$ and hence
$r(e) \geq v_0$.  Thus $| X | > 0$, and by definition $v_0$ is a breaking
vertex.  All that remains is to show that $S = B_H \setminus \{v_0 \}$ . 
Let us suppose that $w_0 \in B_H$.  If $w_0 \notin S$, then $T_{w_0}
\subseteq F^0 \setminus H_S = \{ w \in F^0 \, | \, w \geq \alpha \}$.  But
because the $w_i$'s for $i \geq N_{w_0}$ can only reach elements of $H$ and
$T_{w_0}$, the only way to have $w_i \geq \alpha = \alpha' e_1e_2 \ldots$ for
all $i$ is if we have $w_0 = v_0$.  Hence $v_0$ is the only element of $B_H
\setminus S$ and $S = B_H \setminus \{v_0 \}$.  Thus we have established all
of the claims in part~(2).

For the converse let $E^0 \setminus H$ be a maximal tail.  Consider the
following two cases.
\begin{enumerate}
\item[]{\bf Case I:}  $S = B_H$

We shall show that $F^0 \setminus H_S$ is a maximal tail in $F$.  Since $H_S$
is a saturated hereditary subset of $F^0$, the set $F^0 \setminus H_S$
certainly satisfies (b) and (c) in the definition of maximal tail.  We shall
prove that (a) also holds.  Let $w_1,w_2 \in F^0 \setminus H_S$.  If it is
the case that $w_1,w_2 \in E^0$, then we must also have $w_1,w_2 \in E^0
\setminus H$, and hence there exists $z \in E^0 \setminus H$ such that $w_1
\geq z$ and $w_2 \geq z$ in $E$.  But then $z \in F^0 \setminus H_S$ and $w_1
\geq z$ and $w_2 \geq z$ in $F$.

On the other hand, if one of the $w_i$'s is not in $E^0$, then it must be on
an infinite tail $T_{v_0}$.  Because $w_i \notin H_S$ and $S = B_H$, we must
have $w_i \geq z$ for some $z \in E^0 \setminus H$.  Thus we can replace
$w_i$ with $z$ and reduce to the case when $w_i \in E^0$.

Hence $F^0 \setminus H_S$ also satisfies (a) and is a maximal
tail.  Consequently, $I_{H_S}$ is a primitive ideal by
\cite[Proposition~6.1]{bprs}, and $I_{(H,S)}$ is a primitive ideal by
Corollary~\ref{cor-prim}.

\item[]{\bf Case II:} $E^0 \setminus H = \lambda_{v_0}$ for some breaking
vertex $v_0$ and $S = B_H \setminus \{ v_0 \}$.

As in Case I, it suffices to show that $F^0 \setminus H_S$
satisfies (a) in the definition of maximal tail.  To see this, let $w \in F^0
\setminus H_S$.  If $w \in E^0$, then we must have $w \in E^0
\setminus H = \lambda_{v_0}$ and $w \geq v_0$.  If $w \notin E^0$, then $w$
must be on one of the added tails in $F$.  Since $S = B_H \setminus \{v_0 \}$
we must have that $w$ is an element on $T_{v_0} = \{ v_0, v_1, v_2, \ldots
\}$.  In either case we see that $w$ can reach an element of $T_{v_0}$ in
$F$.  Consequently, $F^0 \setminus H_S \geq T_{v_0}$ and $F^0 \setminus H_S$
clearly satisfies (a).
\end{enumerate}
\end{proof}

\begin{defn}
\label{def-phi-E}
Let $E$ be a graph that satisfies Condition~(K).  We define a map $\phi_E :
\Xi_E \rightarrow \Prim C^*(E)$ as follows.  For $\gamma \in \Lambda_E$ let
$H(\gamma) := E^0 \setminus \gamma$ and define $\phi_E(\gamma) :=
I_{(H(\gamma), B_{H(\gamma)})}$.  For $v_0 \in BV(E)$ we define $\phi_E(v_0)
:= I_{(H(\lambda_{v_0}), B_{H(\lambda_{v_0})} \setminus \{ v_0 \} ) }$.  The
previous theorem shows that $\phi_E$ is a bijection.
\end{defn}

We now wish to define a topology on $\Xi_E$ that will
make $\phi_E$ a homeomorphism.  As usual our
strategy will be to translate the problem to a desingularized graph and
make use of the corresponding results in \cite{bprs}.  In
particular, if $E$ is any graph and $F$ is a desingularization of $E$,  
then we have the following picture: 
$$
\xymatrix{
\Xi_E \ar[d]_{\phi_E} \ar@{.>}[r]^{h} & \Xi_F \ar[d]^{\phi_F}\\
\Prim C^*(E) \ar[r]^{\psi} & \Prim C^*(F),
}
$$
where $\psi$ is the Rieffel correspondence restricted to the primitive ideal
space.  If we use the topology on $\Xi_F = \Lambda_F$ defined in
\cite[Theorem~6.3]{bprs}, then $\phi_F$ is a homeomorphism.  To define a
topology on $\Xi_E$ that makes $\phi_E$ a homeomorphism we will simply use
the  composition $h := \phi_F^{-1} \circ \psi \circ \phi_E$ to pull the
topology on $\Xi_F$ back to a topology on $\Xi_E$.  We start with a
proposition that describes the map $h$. 

\begin{prop}
\label{prop-riefmaxtail}

Let $E$ be a graph satisfying Condition~(K) and let $F$ be a
desingularization of $E$.  
\begin{enumerate}
\item If $\alpha \in E^\infty \cup
\{\alpha
\in E^*
\,|\, r(\alpha) \hbox{ is a singular vertex}\}$ and $\gamma = \{v \in E^0
\,|\, v \geq
\alpha\} \in \Lambda_E$, then
$h(\gamma) = \{v \in F^0 \,|\, v \geq \phi_{\infty}(\alpha)\}$.  
\item If $v_0$ is a breaking vertex, then $h(v_0) = \{ v \in F^0 \, | \,
v \geq e_1e_2 \ldots \}$, where $e_1e_2 \ldots$ is the path on the tail added
to $v_0$.
\end{enumerate}
\end{prop}

\begin{proof}
To prove part (1), let $H := E^0 \setminus \gamma$ and $S := B_H$.  Then using
Proposition~\ref{prop-HSHS} we have $h(\gamma) =
\phi_F^{-1} \circ \psi \circ \phi_E (\gamma) =
\phi_F^{-1} \circ \psi (I_{(H,S)}) = \phi_F^{-1}(I_{H_S}) = F^0 \setminus
H_S$.  We shall show that $F^0 \setminus H_S = \{v \in F^0 \,|\, v \geq
\phi_{\infty}(\alpha)\}$.  To begin, if $v \in E^0$ then 
$$
v \in F^0 \setminus H_S 
\Longleftrightarrow
v \in E^0 \setminus H
\Longleftrightarrow
v \in \gamma
\Longleftrightarrow
v \geq \alpha \text{ in } E
\Longleftrightarrow
v \geq \phi_{\infty}(\alpha) \text{ in } F.
$$
where the last step follows from Lemma~\ref{lem-correspondence}. 
On the other hand, suppose $v \in F^0 \setminus E^0$.  Then since $S=B_H$
every vertex $v \in F^0 \setminus H_S$ must connect to some vertex $w \in E^0
\setminus H$.  So we may replace $v$ with $w$ and repeat the above
argument.  Thus we have proven (1).

For part (2), let $v_0$ be a breaking vertex and set $\lambda_{v_0} := \{ w
\in E^0 \, | \, w \geq v_0 \}$ and $S := B_{(E^0 \setminus \lambda_{v_0})}
\setminus \{ v_0 \}$.  Then $h(v_0) = \phi_F^{-1} \circ \psi \circ \phi_E
(v_0) = \phi_F^{-1} \circ \psi (I_{(H,S)}) = \phi_F^{-1}(I_{H_S}) = F^0
\setminus H_S$.  An argument similar to the one above shows that $F^0
\setminus H_S = \{ v \in F^0 \, | \, v \geq e_1e_2\ldots \}$.
\end{proof}

\begin{defn}
\label{defn-reach}
Let $E$ be a graph and let $S \subseteq E^0$.  If $\gamma$ is a maximal
tail, then we write $\gamma \rightarrow S$ if $\gamma \geq S$.  If $v_0$
is a breaking vertex in $E$, then we write $v_0 \rightarrow S$ if the set
$\{e \in E^0 \,|\, s(e) = v_0, r(e) \geq S\}$ contains infinitely many
elements.
\end{defn}

\begin{lem}
\label{lem-reach}
Let $\delta \in \Xi_E$ and let $P \subseteq \Xi_E$.  Then 
$\delta \rightarrow \bigcup_{\lambda \in P} \lambda$ in $E$
if and only if $h(\delta) \geq \bigcup_{\lambda \in P} h(\lambda)$ in $F$.
\end{lem}

\begin{proof}
If $\delta$ is a maximal tail, then from Lemma~\ref{tailpath} we have $\delta
= \{ v \in E^0 \, | \, v \geq \alpha \}$ for some $\alpha
\in  E^\infty \cup \{\alpha \in E^* \,|\, r(\alpha) \hbox{ is a
singular vertex}\}$.  Similarly, for each $\lambda \in P \cap \Lambda_E$ we may
write $\lambda = \{ v \in E^0 \, | \, v \geq \alpha^\lambda \}$ for some
$\alpha^\lambda \in  E^\infty \cup \{\alpha \in E^* \,|\, r(\alpha) \hbox{ is
a singular vertex}\}$.  Now
\begin{align*}
& \delta \rightarrow \bigcup_{\lambda \in P} \lambda \\
\Longleftrightarrow&
\alpha \geq \bigcup_{\lambda \in P \cap \Lambda_E} \{ r(\alpha^\lambda_i)
\}_{i=1}^{|\alpha^\lambda |} \cup
\bigcup_{v_0 \in P \cap BV(E)} v_0 \\
\Longleftrightarrow&
\phi_\infty(\alpha) \geq \bigcup_{\lambda \in P \cap \Lambda_E}
\{ r(\phi_\infty(\alpha^\lambda)_i) \}_{i=1}^{|\alpha^\lambda|}  \cup
\bigcup_{v_0 \in P \cap BV(E)} \phi_\infty(v_0) \\
\Longleftrightarrow&
\{ v \in F^0 \, | \, v \geq \phi_\infty(\alpha) \} \geq \bigcup_{\lambda \in P
\cap \Lambda_E} \{ v \in F^0 \, | \, v \geq \phi_\infty(\alpha^\lambda_i) \}
\cup \bigcup_{v_0 \in P \cap BV(E)} \{ v \in F^0 \, | \, v \geq
e_1^{v_0}e_2^{v_0}\ldots \} \\ 
\Longleftrightarrow&
h(\delta) \geq \bigcup_{\lambda \in P} h (\lambda)
\end{align*}

\noindent So the claim holds when $\delta$ is a maximal tail.

Now let us consider the case when $\delta = v_0$ is a breaking vertex.  It
follows from Lemma~\ref{prop-riefmaxtail} that $h(v_0) = \{ v \in F^0 \,
| \, v \geq e_1e_2\ldots \}$, where $e_1e_2\ldots$ is the path on the tail
added to $v_0$.  Now suppose that $v_0 \rightarrow \bigcup_{\lambda \in P}
\lambda$.   Fix $v \in h(\delta)$.    Note that either $v \geq v_0$ in $F$ or
$v$ is on the infinite tail added to $v_0$ in $F$.  Because 
$v_0 \rightarrow \bigcup_{\lambda \in P} \lambda$, there are
infinitely many edges in $E$ from $v_0$ to vertices that connect to
$\bigcup_{\lambda \in P} \lambda$.  Thus no matter how far out on the
tail $v$ happens to be, there must be an edge in $F$ whose source is a vertex
further out on the tail than $v$ and whose range is a vertex that connects
to a vertex $w \in \lambda$ for some $\lambda \in P$.  Since $w \in
\lambda$ we must have $w \in h(\lambda)$ and thus $v \geq \bigcup_{\lambda \in
P} h(\lambda)$.

Now assume that 
$h(v_0) \geq \bigcup_{\lambda \in P} h(\lambda)$.  Then
every vertex on the infinite tail attached to $v_0$ connects to a
vertex in $\bigcup_{\lambda \in P} h(\lambda)$.  
In fact it is true that 
every vertex on the infinite tail attached to $v_0$ connects to a
vertex in $\bigcup_{\lambda \in P} h(\lambda) \cap E^0$, which implies
that every vertex on the infinite tail connects to a vertex in
$\bigcup_{\lambda \in P} \lambda$.  But this implies that there must be
infinitely many edges from $v_0$ to vertices that connect to
$\bigcup_{\lambda \in P} \lambda$.  
Thus $v_0 \rightarrow \bigcup_{\lambda \in P} \lambda$.  
\end{proof}

\begin{thm}
\label{thm-topology}

Let $E$ be a graph satisfying Condition~(K).  Then there is a
topology on $\Xi_E$ such that for $S \subseteq \Xi_E$, 
$$
\overline{S} := \{\delta \in \Xi_E \,|\, \delta \rightarrow
\bigcup_{\lambda \in S} \lambda\},
$$
and the map $\phi_E$ given in Definition~\ref{def-phi-E} is a
homeomorphism from $\Xi_E$ onto $\Prim C^*(E)$.

\end{thm}

\begin{proof}
Since $h$ is a bijection, we may use $h$ to pull the topology defined on
$\Xi_F = \Lambda_F$ in \cite[Theorem 6.3]{bprs} back to a topology on
$\Xi_E$.  Specifically, if $S \subseteq \Xi_E$ then
$S=h^{-1}(P)$ for some $P \subseteq \Xi_F$ and we define
$\overline{S} := h^{-1}(\overline{P})$.  But from Lemma \ref{lem-reach}
we see that this is equivalent to defining $\overline{S} = \{\delta
\in \Xi_E \,|\, \delta \rightarrow \bigcup_{\lambda \in S}
\lambda\}$.   Now with this topology $h$, and consequently $\phi_E$, is a
homeomorphism.
\end{proof}

%------------------------------------------------------------------------
\section{Concluding Remarks}
\label{sec-conc}
%-----------------------------------------------------------------------

When we defined a desingularization of a graph in Section~\ref{sec-desing},
for each singular vertex $v_0$ we chose an ordering of the edges
$s^{-1}(v_0)$ and then redistributed these edges along the added tail in such
a way that every vertex on the tail was the source of exactly one of
these edges.  Another way we could have defined a desingularization would be
to instead redistribute a finite number of edges to each vertex on the
added tail.  Thus if $v_0$ is a singular vertex, we could choose a
partition of $s^{-1}(v_0)$ into a countable collection $S_0^{v_0} , S_1^{v_0}
, S_2^{v_0} , \ldots$ of finite (or empty) disjoint sets.  Having done this,
we add a tail to $E$ by first adding a graph of the form 
\begin{equation*}
\xymatrix{
v_0 \ar[r]^{e_1} & v_1 \ar[r]^{e_2} & v_2 \ar[r]^{e_3} & v_3 \ar[r]^{e_4} &
\cdots\\
}
\end{equation*}
We then remove the edges in $s^{-1}(v_0)$ and for each $i$ and each $g
\in S_i^{v_0}$ we draw an edge from $v_i$ to $r(g)$.  More formally, if the
elements of $S_i^{v_0}$ are listed as $\{g_i^1, g_i^2, \ldots, g_i^{m_i} \}$
we define $F^0 := E^0 \cup \{v_i \}_{i=1}^\infty$, $F^1 := (E^1 \setminus
s^{-1}(v_0)) \cup \{e_i \}_{i=1}^\infty \cup \{f_i^j \, | \, 1 \leq i \leq
\infty \text{ and } 1 \leq j \leq m_i \}$, and extend $r$ and $s$ by
$s(e_i) = v_{i-1}$, $r(e_i) = v_i$, $s(f_i^j) = v_i$, and $r(f_i^j) =
r(g_i^j)$.

If we add tails in this manner, then we can define a desingularization of
$E$ to be the graph $F$ formed by adding a tail to each singular vertex in
$E$.  Here a choice of partition $S_0^{v_0} , S_1^{v_0} , S_2^{v_0}
, \ldots$ must be made for each singular vertex, and different choices will
sometimes produce nonisomorphic graphs.

With this slightly more general definition of desingularization, all of the
results of this paper still hold and the proofs of those results remain
essentially the same.  We avoided using this broader definition only because
the partitioning and the use of double subscripts in the $f_i^j$'s creates
very cumbersome notation, and we were afraid that this would obscure the main
points of this article.  However, we conclude by mentioning this
more general method of desingularization because we believe that in practice
there may be situations in which it is convenient to use.  For example, if
$H$ is a saturated hereditary subset of $E^0$, then for each $v_0
\in B_H$ one may wish to choose a partition of $s^{-1}(v)$ with $S_0^{v_0} :=
\{ e \in E^1 \, | \, s(e) = v_0 \text{ and } r(e) \notin H \}$.  Then a
desingularization created using this partition will have the property that
every vertex on a tail added to $v_0$ will point only to the next vertex on
the tail and elements of $H$.

%------------------------------------------------------------------------
%  bibliography
%-----------------------------------------------------------------------
\bibliographystyle{ams}

\begin{thebibliography}{10}

\bibitem{bhrs}
T.~Bates, J.~H.~Hong, I.~Raeburn and W.~Szyma\'nski,
\emph{The ideal structure of $C^*$-algebras of infinite graphs},
preprint (2001).


\bibitem{bprs}
T.~Bates, D.~Pask, I.~Raeburn and W.~Szyma\'nski,
\emph{The $C^*$-algebras of row-finite graphs}, New York J.
Math. \textbf{6} (2000), 307--324.

\bibitem{blackadar}
B.~Blackadar, \emph{Shape theory for ${C}^*$-algebras}, Math. Scand.
  \textbf{56} (1985), 249--275.

\bibitem{Cun}
J.~Cuntz, \emph{Simple $C^*$-algebras generated by
isometries}, Comm. Math. Phys. \textbf{57} (1977), 173--185.

\bibitem{CK1}
J.~Cuntz and W.~Krieger, \emph{A class of ${C}^*$-algebras and topological
  {M}arkov chains}, Invent. Math. \textbf{56} (1980), 251--268.

\bibitem{effros}
E.~Effros, \emph{Dimensions and ${C}^*$-algebras}, CBMS Regional Conf. Ser. in
  Math. no. 46, American Mathematical Society, Providence, RI, 1980.

\bibitem{el}
R.~Exel and M.~Laca, \emph{{C}untz-{K}rieger algebras for infinite matrices},
  J. Reine Angew. Math. \textbf{512} (1999), 119--172.

\bibitem{flr}
N.~Fowler, M.~Laca, and I.~Raeburn, \emph{The ${C}^*$-algebras of infinite
  graphs}, Proc. Amer. Math. Soc. \textbf{8} (2000), 2319--2327.

\bibitem{fr}
N.~Fowler and I.~Raeburn, \emph{The {T}oeplitz algebra of a {H}ilbert
  bimodule}, Indiana Univ. Math. J. \textbf{48} (1999), 155--181.

\bibitem{FW}
M.~Fujii and Y.~Watatani, \emph{Cuntz-Krieger algebras
associated with adjoint graphs},  Math. Japon. \textbf{25},
(1980), 501--506. 

\bibitem{aHr}
A.~an~Huef and I.~Raeburn, \emph{The ideal structure of {C}untz-{K}rieger
  algebras}, Ergo. Th. \& Dynam. Sys. \textbf{17} (1997), 611--624.

\bibitem{kpr}
A.~Kumjian, D.~Pask, and I.~Raeburn, \emph{{C}untz-{K}rieger algebras of
  directed graphs}, Pacific J. Math \textbf{184} (1998), 161--174.

\bibitem{kprr}
A.~Kumjian, D.~Pask, I.~Raeburn, and J.~Renault,
\emph{Graphs, groupoids, and Cuntz-Krieger algebras}, J. Funct.
Anal. \textbf{144} (1997), 505--541.

\bibitem{patgg}
A.~Paterson, \emph{Graph inverse semigroups, groupoids and their
  ${C}^*$-algebras}, J. Operator Theory, to appear.

\bibitem{pimsner}
M.V. Pimsner, \emph{A class of ${C}^*$-algebras generalizing both
  {C}untz-{K}rieger algebras and crossed products by {$\mathbb{Z}$}}, Fields
  Institute Communications \textbf{12} (1997), 189--212.

\bibitem{tfb}
I.~Raeburn and D.P. Williams, \emph{Morita equivalence and continuous trace
  ${C}^*$-algebras}, Math. Surveys and Monographs, no.~60, Amer. Math. Soc.,
  Providence, 1998.

\bibitem{Seg}
I.~Segal, \emph{Irreducible representations of operator
algebras}, Bull. Amer. Math. Soc. \textbf{53} (1947), 73--88.

\bibitem{sz}
W.~Szyma\'nski, \emph{Simplicity of {C}untz-{K}rieger algebras of infinite
  matrices}, Pacific J. Math., \textbf{199} (2001), 249--256. 

\bibitem{Wat}
Y.~Watatani, \emph{Graph theory for $C^*$-algebras}, in Operator
Algebras and Their Applications (R.V. Kadison, ed.), Prpc.
Symp. Pure Math., vol. 38, part~1, Amer. Math. Soc.,
Providence, 1982, 195--197. 


\end{thebibliography}

\providecommand{\bysame}{\leavevmode\hbox to3em{\hrulefill}\thinspace}

\vspace*{.1in}

\end{document}